\documentclass[5p]{elsarticle}
\usepackage{graphicx} 
\usepackage{amsmath}
\usepackage{amssymb}
\usepackage{float}
\usepackage{algorithm}
\usepackage{algpseudocode}
\usepackage{todonotes}

\usepackage{booktabs} 
\usepackage{siunitx} 
\usepackage{multirow}
\usepackage{stfloats}
\usepackage{dsfont}
\usepackage{microtype}
\usepackage[detect-weight=true, detect-family=true]{siunitx} 

\allowdisplaybreaks

\begin{document}

\begin{frontmatter}
    \title{LOTUS: A Warm-Start Framework for Powering Dual Decomposition \\in Large-Scale Two-Stage Stochastic Programs}
    

    \author[a1,a2]{Emma Cornielje}
    \ead{e.j.cornielje@student.vu.nl}

    \author[a1,a2]{Berend Markhorst\corref{c1}}
    \ead{berend.markhorst@cwi.nl}
    

    \author[a2]{Alessandro Zocca}
    \ead{a.zocca@vu.nl}

    \author[a1,a2]{Rob van der Mei}
    \ead{mei@cwi.nl}

    \address[a1]{Stochastics Group, CWI, Amsterdam, The Netherlands}
    \address[a2]{Department of Mathematics, Vrije Universiteit Amsterdam, Amsterdam, The Netherlands}
    
    \cortext[c1]{Corresponding author}
  
    \begin{abstract}
    
    Solving large two-stage stochastic mixed-integer programs is computationally challenging. We propose LOTUS, a subset-based warm-start framework that enhances Dual Decomposition under fixed time budgets. By initializing the dual search with informed multipliers, LOTUS accelerates primal convergence and partially alleviates the impact of weak LP relaxations. Through an extensive computational study on production planning instances, we show that, within two hours, LOTUS yields significantly better primal solutions in 45.83\% of cases, while being outperformed by Dual Decomposition in only 4.17\%.
    
\end{abstract}

\end{frontmatter}

\section{Introduction} \label{sec:introduction}
Stochastic Programming (SP) has become a standard framework for decision making under uncertainty in domains ranging from energy to supply-chain management \citep{applications}. However, accurate representation of uncertainty often necessitates large scenario trees, which makes large-scale two-stage Stochastic Mixed-Integer Programs (SMIPs) computationally intractable due to the curse of dimensionality and non-convexity \citep{a15040103}. While Benders' decomposition effectively handles continuous recourse, SMIPs often require scenario decomposition methods like Dual Decomposition (DD) \citep{DualDecomposition} or Progressive Hedging \citep{PH}.

Despite their usefulness, these methods struggle with convergence at scale. Recent literature suggests that warm-start procedures, which generate initial solutions from reduced problem subsets, can mitigate these computational hurdles \citep{markhorst2024twostepwarmstartmethod}.

In this work, we propose LOTUS (Lagrangian Optimization for Two-stage stochastic programming Using Subsets). LOTUS integrates a warm-start strategy into the DD framework to enhance the discovery of high-quality primal solutions under tight computational budgets. Specifically, we employ Subgradient Descent \citep{subgradient} on a representative subset of scenarios to generate informed dual multipliers, which serve as an effective starting point for the full-scale decomposition. To rigorously benchmark LOTUS and accelerate future research, we introduce a new set of large-scale production planning instances designed to challenge state-of-the-art SMIP solvers.

The remainder of this paper is structured as follows: Section~\ref{sec:literature} reviews related work, Section~\ref{sec:methodology} details the LOTUS framework, and Section~\ref{sec:results} presents computational results demonstrating the framework's attainment of superior primal objective values compared to standard DD within practical time limits.
\section{Literature} \label{sec:literature}
Two-stage stochastic mixed-integer programs are notoriously difficult to solve due to the loss of convexity in the recourse function and the curse of dimensionality associated with the scenario set. While the deterministic equivalent program (DEP) can be solved directly for small instances, large-scale applications require decomposition techniques. These are broadly categorized into stage-wise and scenario-wise decomposition.

\subsection{Decomposition Frameworks} Stage-wise decomposition, and primarily the L-shaped method \citep{OriginalBendersDecomposition,rahmaniani2017benders}, performs well for continuous recourse. However, when integer variables appear in the second stage, the recourse function becomes non-convex and discontinuous. The Integer L-shaped method \citep{LAPORTE1993133} addresses this using specialized cuts (e.g., optimality cuts), but often suffers from slow convergence due to the weakness of these cuts and the computational cost of solving mixed-integer subproblems \citep{af86edcf-19fb-328f-8fc0-79901b485187}.

Scenario-wise decomposition relaxes the non-anticipativity constraints (NACs), thereby splitting the problem into independent scenario subproblems. This structure is well-suited for parallel computation. Two prominent methods in this class are Progressive Hedging (PH) and Dual Decomposition (DD). PH \citep{ProgressiveHedging} applies an augmented Lagrangian approach with a quadratic penalty. While PH is a powerful heuristic for SMIPs, it loses convergence guarantees for integer recourse and does not inherently provide valid dual bounds \citep{LB-PH}.

Conversely, DD \citep{CAROE199937} relaxes the NACs using Lagrangian relaxation without the quadratic penalty. The resulting Lagrangian dual problem provides a valid lower bound (for minimization problems), typically solved via subgradient methods \citep{subgradient} or bundle methods \citep{BundleMethod}. However, for SMIPs, strong duality is not guaranteed; the optimal value of the Lagrangian dual generally provides a lower bound on the primal optimum, often resulting in a duality gap due to the non-convexity of the integer feasible region. To bridge this gap and recover a primal solution, DD is often embedded within a Branch-and-Bound framework (Lagrangian B\&B) \citep{primal}. Despite its theoretical soundness, standard DD is often hindered by the slow convergence of the dual multipliers.

\subsection{Accelerating Convergence: Warm-Starts and Reduction} To address the computational burden of exact methods, recent literature has focused on warm-start strategies \citep{InteriorPoint,markhorst2024twostepwarmstartmethod}. The idea is to use a simplified version of the problem, either by relaxation or by considering a subset of scenarios, to generate advanced initial information for the full-scale problem.

Several frameworks have successfully employed this logic. For interior point methods, \cite{InteriorPoint} demonstrated that solutions from reduced scenario trees provide high-quality starting points. In the context of Branch-and-Cut (B\&C), the TULIP framework \citep{markhorst2024twostepwarmstartmethod} utilizes a two-step approach: solving a problem induced by a subset of scenarios to generate valid inequalities, which then warm-starts the full problem. LOTUS leverages the multi-step framework, but it can be applied to a wider range of two-stage optimization problems. Additionally, \cite{GUO2015311} proposed warm-starting DD with PH, utilizing the rapid heuristic progress of PH to initialize the Lagrange multipliers for the exact DD algorithm.

Our proposed framework, LOTUS, leverages scenario reduction to generate a warm start. Scenario reduction aims to approximate the original probability measure with a discrete measure of smaller support. By approximating the full uncertainty set with a representative subset, LOTUS focuses the search on the core trade-offs of the problem. This filters out the minor scenario variations that typically cause the Master Problem to oscillate, allowing the algorithm to bypass the slow initial search and immediately identify high-quality first-stage decisions for high-quality primal recovery. Scenario reduction methods are classified as either \textit{distribution-based} or \textit{problem-based}. Distribution-based methods, such as fast forward selection \citep{fastforward}, probability distance minimization \citep{probmetrics}, or moment matching \citep{Beltran-Royo02012022}, select scenarios to minimize statistical distance (e.g., Wasserstein metric) without considering the optimization structure of the problem at hand. While computationally fast, they may discard low-probability scenarios that have a high impact on the objective function.

In contrast, problem-driven approaches, see \cite{chou2023problem} for an overview, consider the structure of the underlying optimization problem in their reductions. This idea is used in various algorithms, for example \citep{pdClustering,fsrc,narum2024problem}. LOTUS uses a reduction technique that is both distribution-based and problem-driven to warm-start the DD process, thereby improving the attainment of high-quality primal solutions within practical time limits.
\section{Methodology} \label{sec:methodology}
In this section, we introduce the LOTUS framework. First of all, in Subsection \ref{subsec:method}, each step of the LOTUS algorithm is explained. In Subsection \ref{subsec:assum} we explain the assumptions and motivations for the LOTUS framework. 

Before we explain our methodology, we briefly describe the Deterministic Equivalent Problem (DEP) formulation of a two-stage SMIP. 
\begin{subequations} \label{model:2stage_dep}
\begin{align}
    \min_{x, \{y_{(s)}\}_{s \in \mathcal{S}}} \quad & c^T x + \sum_{s \in \mathcal{S}} p_{(s)} q_{(s)}^T y_{(s)} \label{eq:obj} \\
    \text{s.t.} \quad & A x \geq b, \label{eq:first_stage_con} \\
    & T_{(s)} x + W_{(s)} y_{(s)} \geq h_{(s)}, \quad & \forall s \in \mathcal{S}, \label{eq:second_stage_con} \\
    & x \in \mathbb{Z}_+^{n_1} \times \mathbb{R}_+^{n_2}, \label{eq:x_domain} \\
    & y_{(s)} \in \mathbb{Z}_+^{k_1} \times \mathbb{R}_+^{k_2}, \quad & \forall s \in \mathcal{S}. \label{eq:y_domain}
\end{align}
\end{subequations}
Here $\mathcal{S}$ denotes the set of scenarios, where each scenario $s \in \mathcal{S}$ occurs with probability $p_{(s)}$, with $\sum_{s \in \mathcal{S}} p_{(s)} = 1$. The first-stage decision variables vector is $x$, with cost vector $c$. The first-stage constraints are defined by the matrix $A$ and right-hand side $b$. For each scenario $s$, $y_{(s)}$ denotes the second-stage recourse variables with cost vector $q_{(s)}$. The second-stage constraints are defined by the recourse matrix $W_{(s)}$, the technology matrix $T_{(s)}$, and the right-hand side $h_{(s)}$. Note that these three need not be stochastic. The decision domains defined in \eqref{eq:x_domain} and \eqref{eq:y_domain} represent the mixed-integer requirements, where the total numbers of variables are $n = n_1 + n_2$ and $k = k_1 + k_2$, respectively.

\subsection{Method}\label{subsec:method}
\begin{figure}[t]
    \centering
    \includegraphics[width=1.0\linewidth]{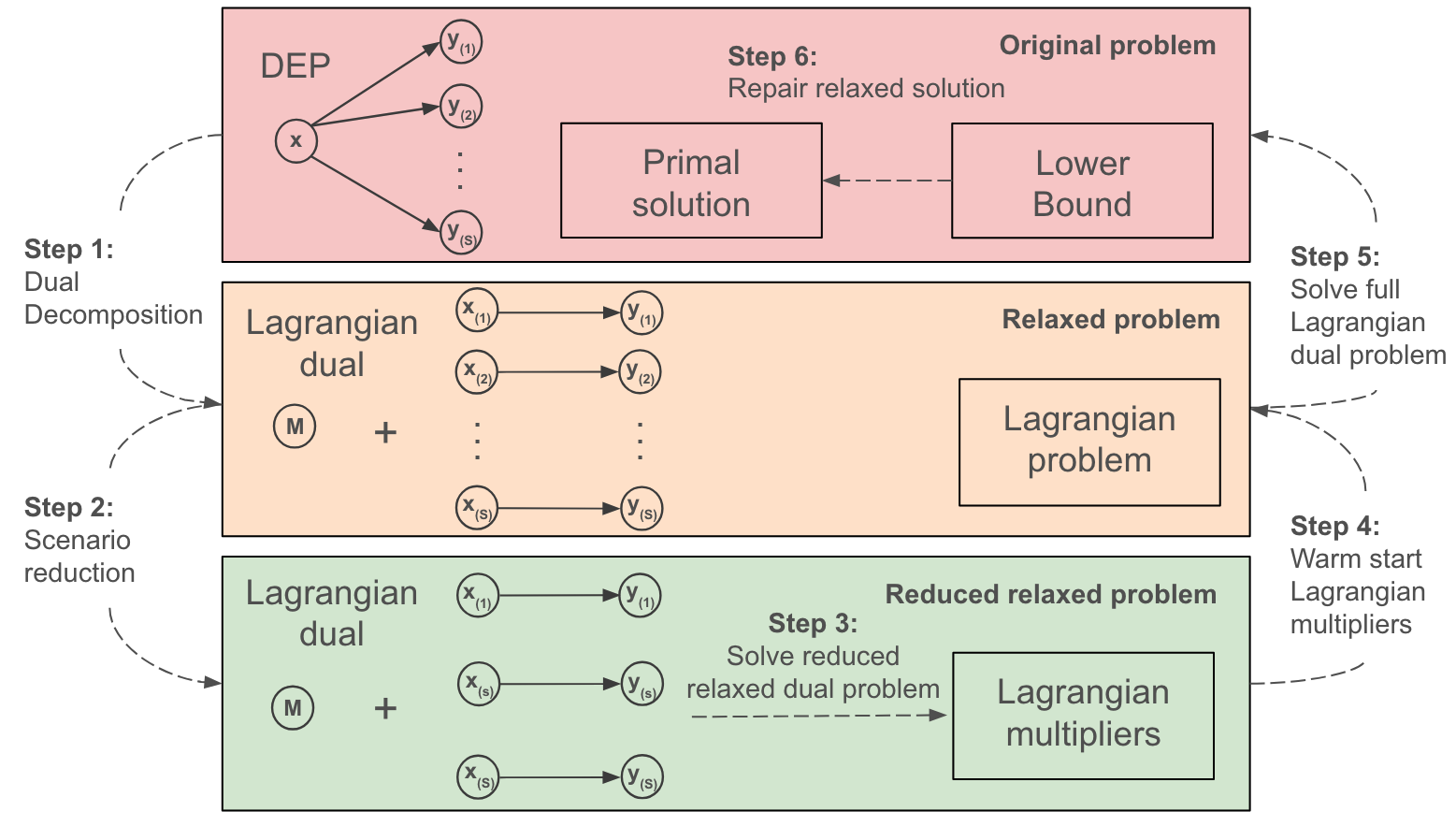}
    \caption{Schematic representation of LOTUS. The dashed arrows correspond to steps in this method. The red, orange, and green boxes denote the original, the full relaxed and the reduced relaxed problem, respectively.}
    \label{fig:LOTUS}
\end{figure}

A schematic overview of the LOTUS framework is given in Figure \ref{fig:LOTUS}. In this overview, the red box corresponds to the original problem, the orange box to the Lagrangian-relaxed formulation of the full problem, and the green box to the Lagrangian-relaxed formulation of the reduced problem. Furthermore, the dashed arrows indicate the steps taken in the framework. Each step applied in LOTUS is discussed in more detail in the following sections.

\subsubsection{Step 1: Dual Decomposition}
We apply \textit{Dual Decomposition} to the DEP formulation. To decouple the scenarios, we introduce local copies of the first-stage variables, $x_{(s)}$, for each scenario $s \in \mathcal{S}$. We enforce explicit non-anticipativity constraints:
\begin{equation*}
    x_{(s)} = x, \quad \forall s \in \mathcal{S}.
\end{equation*}
By relaxing these NACs with Lagrangian multipliers $\lambda_{(s)}$, we obtain the Lagrangian:
\begin{equation*}
    L(x, x_{(s)}, y_{(s)}, \lambda) = c^T x + \sum_{s \in \mathcal{S}} \left( p_{(s)} q_{(s)}^T y_{(s)} + \lambda_{(s)}^T(x_{(s)} - x) \right).
\end{equation*}
Rearranging the terms and minimizing the Lagrangian over $x$, $x_{(s)}$, and, $y_{(s)}$, we obtain the Lagrangian function which decomposes into a \textit{Master Problem} (M) over the global variable $x$, and $|\mathcal{S}|$ independent \textit{Scenario Subproblems} over $(x_{(s)}, y_{(s)})$:
\begin{subequations}
\begin{align}
    \text{(M):} \quad & \min_{x \in \mathcal{X}} \left( c - \sum_{s \in \mathcal{S}} \lambda_{(s)} \right)^T x, \label{eq:master_problem}\\
    \text{(S)}_s: \quad & \min_{x_{(s)}, y_{(s)}} \left( p_{(s)} q_{(s)}^T y_{(s)} + \lambda_{(s)}^T x_{(s)} \right). \label{eq:sub_problem}
\end{align}
\end{subequations}
The Master Problem \eqref{eq:master_problem} minimizes the net first-stage cost over a compact set $\mathcal{X}$ (representing the first-stage feasibility domain, e.g., $Ax \ge b$), preventing unboundedness. The subproblems \eqref{eq:sub_problem} are solved subject to their respective scenario constraints for each $s \in \mathcal{S}$.

For compactness in the subsequent algorithmic steps, we define the \textit{Lagrangian Function} $Z^{\mathcal{S}'}_D(\lambda)$ for any subset of scenarios $\mathcal{S}' \subset \mathcal{S}$ as the minimum of the Lagrangian restricted to that subset:
\begin{equation} \label{eq:dual_func_def}
    Z^{\mathcal{S}'}_D(\lambda) = \min_{x \in \mathcal{X}, \{x_{(s)}, y_{(s)}\}_{s \in \mathcal{S}'}} L_{\mathcal{S}'}(x, x_{(s)}, y_{(s)}, \lambda).
\end{equation}
Consequently, the full Lagrangian dual problem seeks the tightest lower bound by solving:
\begin{equation*}
    \max_{\lambda \in \mathbb{R}^{|\mathcal{S}|}} Z^{\mathcal{S}}_D(\lambda).
\end{equation*}

\subsubsection{Step 2: Scenario Reduction (Warm Start)}
Scenario reduction is the central innovation of the LOTUS framework. Solving the full Lagrangian dual problem is computationally prohibitive for a large collection $\mathcal{S}$ of scenarios. To accelerate convergence, LOTUS generates a \textit{warm start} for the Lagrangian multipliers vector $\lambda$ by solving the dual problem over a much smaller, representative subset of scenarios $\mathcal{S}' \subset \mathcal{S}$, where $|\mathcal{S}'| \ll |\mathcal{S}|$.

The quality of this warm start depends on how accurately $\mathcal{S}'$ approximates the full uncertainty set. We investigate a combination of two metrics for selecting $\mathcal{S}'$:

\paragraph{Distribution-based metric}
Let $\xi_s$ denote the vector of stochastic parameters for scenario $s$. This standard approach characterizes each scenario solely by its data vector $\xi_s$. The distance between scenarios is defined as the Euclidean distance $\|\xi_s - \xi_{s'}\|_2$. While computationally efficient, this metric assumes that proximity in the \textit{data space} implies proximity in the \textit{solution space}, which may not hold for mixed-integer problems where small data perturbations can significantly alter the optimal basis.

\paragraph{Problem-based metric}
To address the limitations of purely data-driven metrics, we introduce a problem-based metric that captures the problem's cost structure. In MIPs, small changes in data can trigger large discontinuities in the objective value (e.g., exceeding a capacity threshold that incurs a fixed cost).
Therefore, we map each scenario $s$ to a scalar feature $v_s^{LP}$, representing the optimal objective value of its linear programming relaxation. This is obtained by solving the scenario subproblem individually (ignoring non-anticipativity constraints) with all integer requirements relaxed. The distance between scenarios, e.g., $s$ and $s^{'}$, is then defined based on the difference in these structural proxy values, $|v_s^{LP} - v_{s'}^{LP}|$.

\paragraph{Combining the metrics}
We combine the data vector $\xi_s$ and the feature $v_s^{LP}$ into one feature vector $\mathcal{F}_s$. To ensure statistical balance between distribution-based and problem-based features, we standardize both using $Z$-score normalization before concatenation. 
To ensure the problem-based metric and distribution-based metric carry equal weight in the selection algorithm, we apply a weight $\omega$ to the normalized cost feature equal to the number of furniture types ($|F|$). This balancing ensures that the selection algorithm considers the underlying cost impact of scenarios with the same total priority, given the multi-dimensional demand data.

\paragraph{Selection algorithm}
We then employ a Fast Forward Selection algorithm to iteratively build the subset $\mathcal{S}'$ by selecting the scenario that minimizes the probabilistic distance (e.g., Kantorovich distance) between the original distribution and the reduced distribution. Once $\mathcal{S}'$ is fixed, the original set $\mathcal{S}$ is re-partitioned: every discarded scenario $k \in \mathcal{S} \setminus \mathcal{S}'$ is mapped to its closest representative $s \in \mathcal{S}'$ (according to the Euclidean distance metric $\|\mathcal{F}_s - \mathcal{F}_{s'}\|_2$), and their probabilities are aggregated to preserve the total probability mass.

\subsubsection{Step 3: Solve the Reduced Lagrangian Dual}
In this step, we solve the Lagrangian dual problem restricted to the representative subset of scenarios $\mathcal{S}' \subset \mathcal{S}$:
\begin{equation*} 
    \max_{\lambda \in \mathbb{R}^{|\mathcal{S}|}} Z_D^{\mathcal{S}'}(\lambda).
\end{equation*}
Using the notation defined in Step 1, this seeks the tightest lower bound for the reduced problem. Since $Z_D^{\mathcal{S}'}(\lambda)$ is concave and non-differentiable, we employ the \textit{Subgradient Ascent algorithm} \citep{subgradient}. For further details, see~\ref{sec:details_step_3}.

\subsubsection{Step 4: Warm Start Initialization}
The optimized multipliers $\lambda^*$ obtained from the reduced problem in Step 3 are mapped to the full scenario set $\mathcal{S}$ to initialize the full-scale dual problem. Let $\rho: \mathcal{S} \to \mathcal{S}'$ be the mapping function from the scenario reduction step, where $\rho(i)$ denotes the representative scenario in $\mathcal{S}'$ for any original scenario $i \in \mathcal{S}$. The initial multipliers $\lambda^0$ for the full problem are set as: $\lambda^0_{(i)} = \lambda^*_{(\rho(i))}$ for every $i \in \mathcal{S}$. This warm-start configuration provides the initial point for the Subgradient Ascent algorithm applied to the full Lagrangian dual problem.

\subsubsection{Step 5: Solve the Full Lagrangian Dual}
We now solve the full-scale Lagrangian dual problem over the complete scenario set $\mathcal{S}$:
\begin{equation*}
    \max_{\lambda} Z_D^{\mathcal{S}}(\lambda).
\end{equation*}

\paragraph{Initialization (warm start)}
Unlike in Step 3, we use the information obtained from the reduced problem. We initialize the multipliers using the \textit{warm start} vector $\lambda^{\text{WS}}$ constructed in Step 4, i.e.,
$\lambda^0 = \lambda^{\text{WS}}$. This advanced starting point in the dual space is expected to be significantly closer to the optimal multipliers than the zero vector, thereby improving primal convergence by identifying superior feasible solutions earlier in the optimization process.

\paragraph{Algorithm}
The remaining algorithmic procedure (Subgradient Ascent, Polyak step size, and termination criteria) is identical to that described in Step 3, ensuring a consistent approach to solving the full dual problem. Let $\lambda^*$ denote the optimal (or best found) multipliers returned by the algorithm upon termination.

\subsubsection{Step 6: Primal Recovery and Evaluation} \label{subsubsec:primal}
While the Lagrangian dual function $Z_D^{\mathcal{S}}(\lambda^*)$ provides a rigorous lower bound on the optimal value, the resulting scenario-specific solutions $\{x_{(s)}^*(\lambda^*)\}_{s \in \mathcal{S}}$ typically violate non-anticipativity (i.e., $x_{(s)}^* \neq x_{(s')}^*$). To recover a feasible \textit{primal} solution $Z_P^{\mathcal{S}}(\lambda^*)$ (and thus a valid upper bound for a minimization problem), we leverage the global decisions $\tilde{x}$ from the Master Problem to conduct a restricted primal search of the full DEP.

In this step, we solve the DEP with the first-stage variables constrained to a narrow feasibility window $(1 \pm \epsilon)\tilde{x}$ with $\epsilon = 0.05$. This approach uses the progress of the dual search to prune the primal space, focusing the solver on the most promising region of the first-stage domain $\mathcal{X}$. Due to the assumption of relatively complete recourse, see Section~\ref{subsec:assum}, it is ensured that this restricted DEP will find optimal second-stage recourse variables $\{y_{(s)}\}_{s \in \mathcal{S}}$ for every scenario. The resulting objective value provides a valid primal bound $Z_P^{\mathcal{S}}$, which we use to evaluate the final optimality gap.

\paragraph{Note on optimality}
If the gap between the primal bound and the dual bound is non-zero, the solution is not guaranteed to be optimal. To ensure exact convergence for non-convex SMIPs, LOTUS is designed to be embedded as a bounding procedure within, e.g., a global B\&B scheme \citep{primal}. The accelerated primal convergence provided by LOTUS facilitates more aggressive pruning of the search tree early on, thereby reducing the number of nodes explored and mitigating the computational burden of weak dual bound convergence.

\subsection{Theoretical Assumptions and Practical Motivations} \label{subsec:assum}
The design of the LOTUS framework is driven by specific structural properties of a Two-Stage SMIP. In this section, we outline the theoretical conditions for LOTUS and the practical challenges that motivate the proposed warm-start methodology.

\subsubsection{Theoretical Assumptions}
For LOTUS, we make the following assumptions regarding the problem structure.

\paragraph{Block-angular structure}
We assume the constraint matrix has a block-angular structure, in which scenarios are linked solely through the first-stage decision variables $x$. This separability is the fundamental property that enables the Lagrangian relaxation of the NACs to decompose the DEP into independent, tractable subproblems.

\paragraph{Relatively complete recourse and compactness}
We assume the problem exhibits \textit{relatively complete recourse} and that the feasible regions $\mathcal{X}$ and $\mathcal{Y}_s$ are \textit{compact} (closed and bounded). These properties are critical for two reasons:
\begin{enumerate}
    \item \textit{Algorithmic Stability:} The Lagrangian dual function $Z_D^{\mathcal{S}}(\lambda)$ is \textit{well-defined} and finite for all $\lambda \in \mathbb{R}^n$. This prevents infinite penalties in the objective or unbounded solutions from destabilizing the subgradient updates.
    \item \textit{Primal Recovery:} It guarantees that the projection heuristic in Step 6 yields a feasible solution, as fixing $x$ to a feasible $\hat{x}$ will not result in infeasibility in the second stage.
\end{enumerate}

\subsubsection{Practical motivations}
Beyond structural validity, the LOTUS framework is specifically designed to address the limitations of standard solvers and classical decomposition methods for SMIPs.

\paragraph{Tightening the integrality gap}
Standard B\&B solvers rely on linear programming relaxations to prune the search tree. However, SMIPs often suffer from \textit{weak LP relaxations}, where the dual bound (lower bound) $Z_D^{\mathcal{S}}$ provided by the relaxed LP is far from the optimal integer value $Z_I^{\mathcal{S}}$. This large integrality gap leads to ineffective pruning and intractable runtimes. DD is motivated by the fact that the Lagrangian dual bound $Z_D^{\mathcal{S}}$ is strictly tighter than (or at least equal to) the LP relaxation bound $Z_{LP}^{\mathcal{S}}$:
\begin{equation*}
    Z_{LP}^{\mathcal{S}} \le Z_D^{\mathcal{S}} \le Z_I^{\mathcal{S}} \le Z_P^{\mathcal{S}}.
\end{equation*}
By solving subproblems that retain their integrality constraints, LOTUS generates a high-quality dual bound $Z_D^{\mathcal{S}}$. More importantly, the framework leverages stable multiplier guidance to prioritize the discovery of high-quality feasible solutions ($Z_P^{\mathcal{S}}$), effectively narrowing the optimality gap from the top down even when the theoretical dual bound remains relatively loose.

\paragraph{Accelerating primal convergence via warm start}
Standard subgradient methods are notoriously slow to converge, largely due to the ``identification phase'': they typically spend many initial iterations simply identifying the non-zero components of the optimal multiplier vector $\lambda^*$ (i.e., the active manifold).
LOTUS circumvents this bottleneck by performing this identification on the computationally cheaper reduced problem $Z_D^{\mathcal{S}'}(\lambda)$. The resulting warm start $\lambda^{\text{WS}}$ effectively projects the full-scale algorithm past the initial search phase, starting it closer to the neighborhood of the optimal active set. This proximity reduces the erratic oscillations typical of cold-start subgradient methods, providing more stable guidance for the primal heuristics.
\section{Results} \label{sec:results}

\subsection{Experimental setup}
To empirically evaluate the performance of the LOTUS algorithm compared to DD, we conducted a two-phase computational case study in the context of a production/manufacturing problem (see \ref{sec:case_study} for details). 

\subsection{Empirical performance of LOTUS}
We evaluate the performance of the LOTUS framework against a standard DD baseline. The primary goal is to assess whether the warm-start strategy provides a significant advantage in finding high-quality primal solutions and closing the optimality gap within a practical time limit of $120$ minutes. For descriptive metrics, we use a stabilized analysis across $24$ unique configurations, where the result of each configuration is the mean of $5$ independent seeds. For a statistical analysis, we perform the Wilcoxon signed-rank test on the $N = 111$ raw instances \cite{wilcoxon} (based on the ranks of the differences) to determine whether the median difference for the primal bound and the optimality gap between DD and LOTUS is significantly different from zero. Note that $9$ instances were removed from the statistical analysis. These are instances in which DD failed to identify a feasible positive primal solution.

Table \ref{tab:overall_sum} below shows that LOTUS achieved an average profit improvement of $0.93\%$. The $95\%$ Confidence interval confirms the robustness of this improvement, as the entire interval remains strictly above zero. While the Average Profit Improvement may appear modest, the global $p$-value of $<0.001$ confirms that LOTUS consistently identifies superior primal solutions. LOTUS demonstrates reliability, providing a solution equal (draw) to or better (win) than standard DD in $95.8\%$ (win rate + draw rate) of tested configurations. 

\begin{table}[H]
\centering
\caption{Overall performance summary ($N = 111$)}
\label{tab:overall_sum}
\small 
\setlength{\tabcolsep}{1.7pt} 
\begin{tabular}{@{}cccccc@{}}
\toprule
\shortstack{Avg.~Profit\\Improv.~(\%)} & \shortstack{Win\\Rate (\%)} & \shortstack{Draw\\Rate (\%)} 
& \shortstack{95\% Conf.\\Interval} &
\shortstack{Primal\\$p$-val} &
\shortstack{Gap\\$p$-val}
\\ \midrule
$0.93$ & $45.83$ & $50.00$ & $[0.38,1.49]$ & $\mathbf{< 0.001}$ & $0.349$\\ \bottomrule
\end{tabular}
\end{table} 

However, we observe a statistical divergence between the objective value and the optimality gap. The $p$-value for the optimality gap ($0.349$) suggests that LOTUS does not necessarily close the gap more efficiently than DD. This is consistent with our practical motivation: LOTUS prioritizes the top-down narrowing of the gap by improving the primal bound, even when the theoretical dual bound remains relatively stagnant.
\begin{table}[H]
\centering
\caption{Performance split by number of scenarios ($|\mathcal{S}|$) over $2$-hour and $4$-hour runtimes.}
\label{tab:results_ns}
\small 
\setlength{\tabcolsep}{3.5pt} 
\begin{tabular}{@{}lrcrcrrr@{}}
\toprule
$|\mathcal{S}|$ & \shortstack{Avg. Profit \\Improv.} & \shortstack{Win\\Rate} & \shortstack{Draw\\Rate} & $N$ & \shortstack{Primal \\$p$-val} & \shortstack{Gap\\$p$-val} & Time \\ \midrule
$250$ & $1.94\%$ & $83.33\%$ & $0.0\%$  & $28$ & $\mathbf{0.0015}$ & $0.0026$ & 2h \\
$500$ & $0.93\%$ & $66.67\%$ & $33.33\%$ & $29$ & $\mathbf{0.0006}$ & $0.2085$ & 2h \\
$1000$ & $0.45\%$ & $33.33\%$ & $66.67\%$ & $28$ & $\mathbf{0.0213}$ & $0.6772$ & 2h \\
$2000$ & $0.41\%$ & $16.67\%$ & $83.33\%$ & $26$ & $0.0899$ & $0.9898$ & 2h \\
$2000$ & $0.67\%$ & $33.33\%$ & $66.67\%$ & $26$ & $\mathbf{0.0090}$ & $0.3383$ & 4h \\ \bottomrule
\end{tabular}
\end{table}
Table \ref{tab:results_ns} details the performance of LOTUS as a function of the number of scenarios ($|\mathcal{S}|$). LOTUS is most dominant at smaller scales ($|\mathcal{S}| = 250$), achieving a $1.94\%$ increase in profit and winning in nearly $83\%$ of cases. As the problem scale increases to $|\mathcal{S}| = 2000$, the frequency of draws increases to $83.33\%$ within the $2$-hour run limit, and the average profit improvement is not significant anymore. Extending the runtime to $4$ hours for the $|\mathcal{S}| = 2000$ instances, allowing for a more exhaustive warm-start phase ($30$ iterations instead of $10$ iterations), reveals the framework's latent efficacy, yielding a significant $p$-value of $0.009$ and an increase in the win rate to $33.33\%$. This indicates that while a larger instance requires a longer warm-start phase, the superior coordination provided by the LOTUS warm start ultimately leads to improved primal convergence.

A critical finding is that LOTUS achieves a higher total iteration count than DD across all tested numbers of scenarios. At $|\mathcal{S}|=2000$, LOTUS completes $50$ total iterations ($10$ during the warm-start and $40$ during the main phase), whereas DD only performs $36$ iterations. This iteration surplus confirms that LOTUS does not waste time on warm-starting. LOTUS seems to trade expensive full-scale iterations for cheaper reduced-scale iterations.

\begin{figure}[t]
    \centering
    \includegraphics[width=\linewidth]{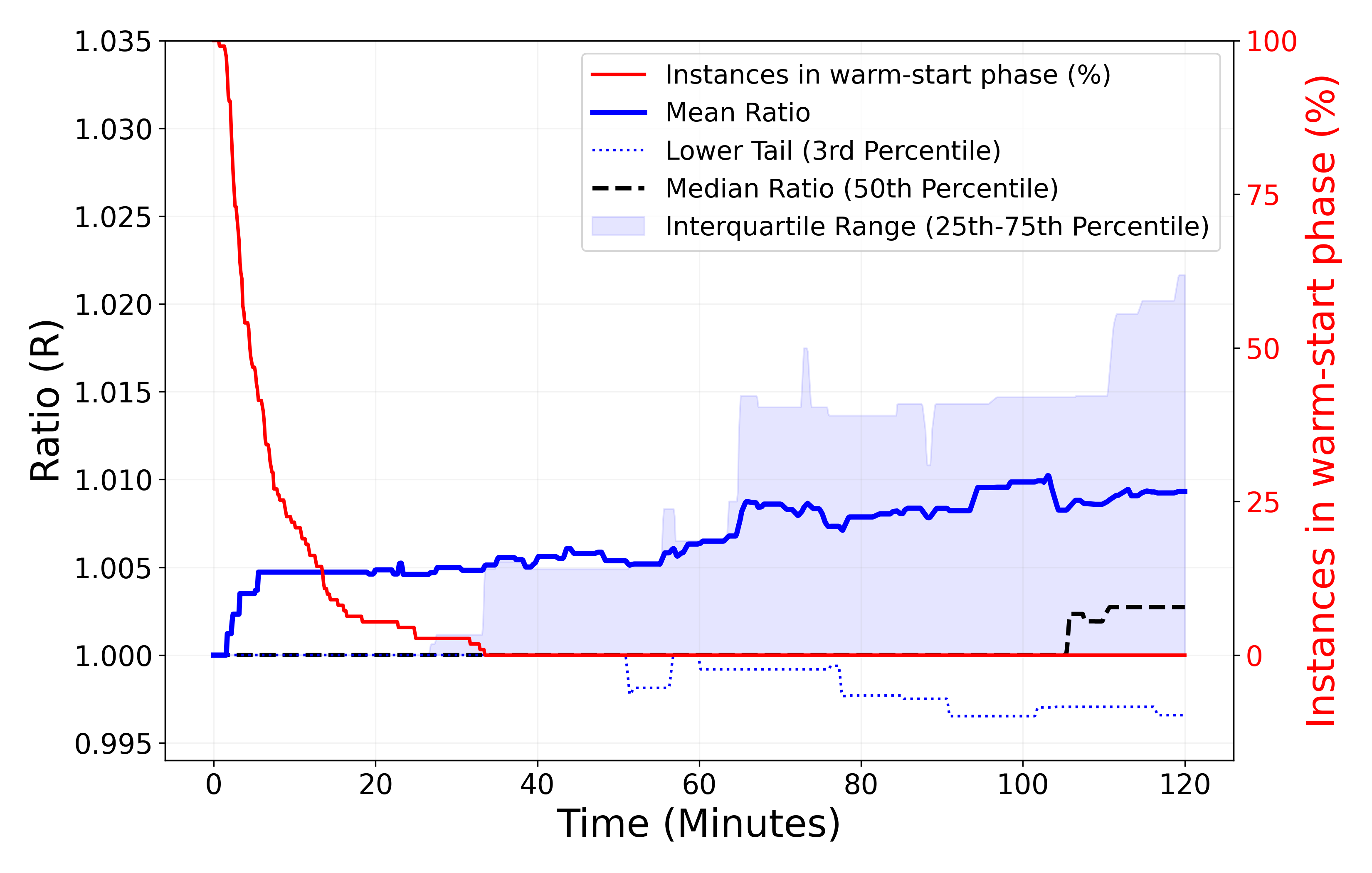}
    \caption{Aggregated primal bound ratio ($R = Z_{P,LOTUS} / Z_{P,DD}$) over time (in blue) and percentage of instances in the warm-start phase over time (in red).}
    \label{fig:lb_ratio}
\end{figure}

To understand the performance of LOTUS relative to DD over the $120$-minute horizon, we analyze the relative primal solution (lower bound) ratios over the full runtime. Figure \ref{fig:lb_ratio} illustrates the evolution of the ratio $R := Z_{P,LOTUS} / Z_{P,DD}$ over time. 
Figure \ref{fig:lb_ratio} shows that the mean ratio $R$ (solid blue line) starts at exactly $1.000$ when $100\%$ of the instances are in the warm-start phase (see solid red line). As instances exit the warm-start phase (as shown by the drop in the red line), they enter the full-scale optimization phase, in which the mean ratio exhibits an upward trajectory. As each instance completes its specific warm-start exploration, it starts the main phase with a superior initial primal bound. By the time the red line reaches $0\%$, at which all configurations have transitioned to solving the full scenario set, LOTUS has secured an average head-start of $0.48\%$ over DD in terms of the mean ratio $R$. This confirms that the warm-start multipliers $\lambda^{WS}$ allow LOTUS to overcome early convergence plateaus where standard DD stagnates. While the frequency of draws ($45.83\%$) keeps the median near $1.000$ for most of the time horizon, the $75$-th percentile diverges significantly after $30$ minutes. Furthermore, the higher median values near the end of the $2$-hour runtime indicate that LOTUS identifies superior solutions during the final convergence phase. However, since the framework incurs a lower profit in only $4.17\%$ of cases, the $25$-th percentile remains at or above $1.000$. This means that the underperformance of LOTUS relative to DD becomes visible only in the extreme lower tail ($3$-th percentile). 

Figure \ref{fig:lb_ratio} confirms that for the vast majority of instances, LOTUS is at least as robust as standard DD, while providing a significant improvement in nearly half of the cases. The results also confirm that LOTUS exploits a double advantage: the warm-start subgradient ascent iterations are much cheaper than in the full problem ($t_L(S') \ll t_L(S)$), and in the main phase, a higher subgradient ascent iteration velocity is achieved. This combined effect generates the primal bound improvement compared to DD seen in Figure \ref{fig:lb_ratio}, resulting in a better primal mean and median for LOTUS within the same $120$-minute window.

\ref{sec:theoretical_insights} presents further insight and empirical analysis of the speed gains of LOTUS with respect to DD.
\section{Discussion} \label{sec:discussion}
While the results for the LOTUS are promising, several limitations must be acknowledged to contextualize the findings and guide future development.

\subsection{Experimental design}
A primary limitation is that warm-start configurations, such as the subset size and number of warm-start iterations, were selected based on preliminary experiments and applied universally. This assumes a degree of scale independence that may not hold for extreme instances. As the problem dimensionality increases, the optimal balance between the warm-start and main decomposition phases shifts. For our largest instances ($|\mathcal{S}|=2000$), the two-hour limit proved insufficient to fully recover the initial computational overhead of the warm-start, leading to a high frequency of draws. However, the emergence of a significant $p$-value ($0.009$) at the four-hour mark confirms that LOTUS has the potential for performance gains, but this depends on the instance size and the allocated warm-start time. The existence of this ``break-even point'', the moment when the time invested in low-cost exploration is recovered through superior primal convergence, suggests that future research is required for dynamic parameter tuning or adaptive subset scaling.

\subsection{Algorithmic framework and benchmarking}
The analysis relies on a standard Subgradient Descent implementation for both the warm-start and the main phase. While this isolates the impact of LOTUS, subgradient methods are known for their sensitivity to step-size rules and potential for slow tail-end convergence. More sophisticated dual-update mechanisms, such as bundle methods \citep{BundleMethod}, could offer better stability and faster dual bound progression. Furthermore, while standard DD was the appropriate benchmark for this foundational study, future work should integrate LOTUS with commercial-grade heuristics to compare its performance directly against the state-of-the-art in industrial solvers.

\subsection{Generalizability and problem structure}
The generalizability of these findings is linked to the structural properties of the set of test instances. The production planning instances include common SMIP challenges, such as fixed-charge costs, integrality in the second stage, and weak LP relaxations. On the other hand, the quantitative gains (the $0.93\%$ average profit improvement) are influenced by the model's specific recourse structure. In problems where the first-stage and second-stage decisions are more loosely coupled, or where the cost of uncertainty is lower, the impact of a high-quality dual warm-start may be less pronounced. 

Testing LOTUS in other domains with similar mathematical problem structures as described in Section \ref{subsec:assum} will be necessary to establish a more universal performance profile. Promising candidates for future experiments include stochastic unit commitment in power systems, which is characterized by weak LP relaxations stemming from significant fixed-charge start-up costs and a reliance on scenario-based decomposition \citep{CAROE199937, papavasiliou2013multiarea}. Similarly, capacitated fixed-charge network design problems exhibit large integrality gaps and dense coupling, which LOTUS is designed to mitigate \citep{gendron2018node}. In both domains, LOTUS's ability to identify a stable active manifold on a reduced set of scenarios could potentially overcome the computational bottlenecks currently faced by standard DD.
\section{Conclusion} \label{sec:conclusion}

Solving large-scale, real-world two-stage SMIPs presents considerable computational challenges stemming from their non-convexity, high dimensionality, and large LP relaxations. In this work, we introduced and evaluated the Lagrangian Optimization for Two-stage stochastic programming Using Subsets (LOTUS) framework, a novel approach that integrates a subset-based warm-start strategy into the Dual Decomposition (DD) method. By initializing the dual search with informed multipliers derived from a representative scenario subset, LOTUS significantly enhances the discovery of high-quality primal solutions within practical computational windows. The reason why this happens is two-fold: first, its warm-start iterations are cheap, and, second, it achieves a higher iteration rate when going back to solving the full problem. 

The results of an extensive computational study across $120$ test instances demonstrate that LOTUS is a robust and practical tool for industrial decision-making. While the mathematical optimality gap remains comparable to standard DD, LOTUS identifies significantly superior primal bounds ($p < 0.001$), achieving an average profit improvement of $0.93\%$. Furthermore, with a $95.8\%$ success rate, the framework offers a level of reliability that DD lacks in the early stages of the solve, where stable guidance is critical for integer recovery.

Our analysis reveals a scale-dependent performance profile. While LOTUS provides immediate benefits for medium-sized problems, the largest instances require longer execution windows to overcome the initial computational overhead of the warm-start phase. The shift from frequent draws at the two-hour mark to a significant win rate at the four-hour mark for the largest instances confirms that the warm start provides superior coordination, eventually enabling the solver to break through the primal convergence plateaus where DD stagnates.

The limitations identified in Section~\ref{sec:discussion} give rise to several promising directions for future research. First, LOTUS could be integrated into a complete solution methodology, such as a B\&B or B\&C framework. The accelerated primal convergence provided by LOTUS facilitates more aggressive pruning of the search tree earlier on, thereby reducing the number of nodes explored and mitigating the computational burden of weak dual bound convergence. Furthermore, replacing subgradient ascent with more stable bundle or level methods could potentially accelerate the warm-start phase. 

Finally, developing an adaptive LOTUS framework in which the warm-start configuration and subset sizes are dynamically determined based on convergence behavior could further optimize the temporal break-even point between the initial warm-start and the main execution phase.\\

{\small \textbf{Acknowledgments}
This publication is partly financed by the Dutch Research Council (NWO) through the ``Readiness'' project (\textit{TWM.BL.019.002}) and the personal VIDI grant of A.~Zocca (\textit{VI.VIDI.233.247}). The authors have no competing interests to declare that are
relevant to the content of this work.}\\

{\small \textbf{Declaration of Use of Generative AI}
The authors used generative artificial intelligence tools solely for the purpose of improving the clarity and linguistic quality of the manuscript. The tools were not used to generate scientific content, results, proofs, data, or interpretations. All substantive content reflects the authors’ original work, and the authors take full responsibility for the accuracy and integrity of the manuscript.}

\bibliographystyle{abbrvnat} 
\small \bibliography{MyReferencesFile,references_berend}

\newpage

\appendix

\section{Details Step 3 (Subgradient Ascent)} \label{sec:details_step_3}
\paragraph{Initialization (cold start)}
Since no prior information is available, we initialize the multipliers with a \textit{cold start}:
$\lambda^0_{(s)} = 0$ for all $s \in \mathcal{S}$.

\paragraph{Step size rule}
To determine the step size $\alpha_k$ at iteration $k$, we utilize the Polyak step size rule \citep{polyak}:
\begin{equation}
    \alpha_k = \gamma_k \frac{Z^{\mathcal{S}'}_{P} - Z^{\mathcal{S}'}_{D}(\lambda^k)}{\|g^k\|^2},
\end{equation}
where $Z_D^{\mathcal{S}'}(\lambda^k)$ is the current dual objective value, $Z_{P}^{\mathcal{S}'}$ is the best known upper bound (primal feasible solution) found so far, and $\gamma_k \in (0, 2]$ is a scaling parameter initialized at $\gamma_0 = 1.8$.

\paragraph{Primal heuristic (trust region)}
To tighten the upper bound $Z_{P}^{\mathcal{S}'}$ needed for the Polyak rule, we run a primal heuristic every 20 iterations. We employ a \textit{fix-and-optimize} strategy: we define a Restricted Master Problem (RMP) by enforcing a 5\% tolerance neighborhood (trust region) around the current best first-stage solution $\bar{x}$. This significantly reduces the search space, allowing a MIP solver to rapidly identify high-quality feasible solutions.

\paragraph{Adaptive tuning and termination}
We implement an adaptive mechanism to dampen oscillations: if the best dual bound fails to improve for five consecutive iterations, we update $\gamma_{k+1} \leftarrow 0.5 \gamma_k$. The algorithm terminates when either a maximum runtime $T_{\text{max}}$ is reached or the relative optimality gap falls below $\epsilon = 0.01\%$.

\section{Case study} \label{sec:case_study}
\subsection*{Model formulation}
To empirically evaluate the performance of LOTUS, we designed a comprehensive computational study utilizing the production planning problem introduced in \cite{Higle2005}. The problem concerns a manufacturer that produces furniture (e.g., desks, tables, chairs) using a finite amount of resources (e.g., lumber, finishing, carpentry) with known demand. We extend this problem to a stochastic setting by treating market demand as an uncertain parameter. We also extend the basic model of \cite{Higle2005} to incorporate several features that reflect real-world operational constraints. While the base model remains tractable for standard solvers, these extensions introduce significant computational challenges. Specifically, the inclusion of binary variables for fixed costs, integrality constraints for production batching, and penalty logic for service failures results in a highly non-convex search space that is difficult for traditional B\&B methods to navigate effectively. We obtain the following DEP formulation:
\begin{small}
\begin{subequations}
\begin{alignat}{3}
    \min_{\bar{x},\{y_{(s)}\}_{s \in \mathcal{S}}, \{\upsilon_{(s)}\}_{s \in \mathcal{S}}, \alpha} \quad & c^T\bar{x}
    - \sum_{s \in \mathcal{S}} p_{(s)} q^T y_{(s)} \label{eq:objB}\\
    & + \sum_{s \in \mathcal{S}} p_{(s)} f^T \upsilon_{(s)} + u^T \bar{\alpha} \notag
\end{alignat}
    \begin{alignat}{3}
    \text{s.t.} \quad  & x_{(s)} \geq Wy_{(s)}, \quad & \forall s \in \mathcal{S} \\
    & y_{(s)} + \upsilon_{(s)} = d_{(s)}, \quad & \forall s \in \mathcal{S}  \label{eq:demand_balance}\\
    & x_{(s)} \leq \bar{x}, \quad & \forall s \in \mathcal{S} \\
    & \alpha_{(s)} \leq \bar{\alpha}, \quad & \forall s \in \mathcal{S} \label{eq:nacalpha}\\
    & \bar{x} \leq L^T \alpha, \label{eq:fixed_cost_link}\\
    & y_{(s)} \geq b \beta_{(s)}, \quad & \forall s \in \mathcal{S} \label{eq:min_batch}\\
    & y_{(s)} \leq d_{(s)} \beta_{(s)}, \quad & \forall s \in \mathcal{S} \label{eq:prod_indicator}\\
    & x_{(s)} \geq 0, \quad & \forall s \in \mathcal{S} \\
    & \bar{x} \geq 0, \\
    & y_{(s)} \in \mathbb{Z}_+^{|F|}, \quad & \forall s \in \mathcal{S} \\
    & \upsilon_{(s)} \in \mathbb{Z}_+^{|F|}, \quad & \forall s \in \mathcal{S} \\
    & \alpha \in \{0,1\}^{|R|}, \\
    & \bar{\alpha} \in \{0,1\}^{|R|}, \\
    & \beta_{(s)} \in \{0,1\}^{|F|}, \quad & \forall s \in \mathcal{S}. 
\end{alignat}
\end{subequations}
\end{small}
We modeled the uncertainty using a finite set of scenarios, $s \in \mathcal{S}$. Each scenario represents a potential future demand state $d_{(s)}$, and is assigned a probability $p_{(s)}$, such that $\sum_{s \in \mathcal{S}} p_{(s)} = 1$. The decision variables in the model are the vector of acquired resources $\bar{x} \in \mathds{R}_+^{|R|}$, the vectors of produced furniture units $y_{s} \in \mathds{Z}_+^{|F|}$, $\forall s \in \mathcal{S}$, and the vectors of unmet demand $\upsilon_{s} \in \mathds{Z}_+^{|F|}$, $\forall s \in \mathcal{S}$. The objective function \eqref{eq:objB} minimizes the net cost over all scenarios, which is the sum over all scenarios of the total variable cost of resources $c^T\bar{x}$ minus the total sales revenue $q^Ty_{s}$ plus the total fixed cost of resources $ u^T \alpha$ plus the
expected penalty for unmet demand $f^T \upsilon_{(s)}$. The model relies on several key parameters. Vectors $c \in \mathds{R}_+^{|R|}$ and $q \in \mathds{R}_+^{|F|}$ contain the costs of acquiring each resource and the sales prices for each type of furniture, respectively. The market limitations are captured by the demand vectors $d_{s} \in \mathds{R}_+^{|F|}, \forall s \in \mathcal{S}$, which specify the maximum number of units that can be sold for each item in each scenario. The production process itself is defined by the technology matrix $W \in\mathds{R}_+^{|R| \times |F|}$, where each element $w_{rf}$ specifies the amount of resources $r$ needed to produce a unit of furniture $f$. The constraint \eqref{eq:demand_balance} defines the service shortage. The constraint \eqref{eq:fixed_cost_link} enforces the limit of resources and links resource acquisition to the fixed cost variable $\alpha$. Constraints \eqref{eq:min_batch} and \eqref{eq:prod_indicator} together manage the minimum batch-size logic for production.

\paragraph{Implementation note} While the general framework in Section \ref{sec:methodology} employs equality NACs ($x_{(s)} = x$), we specialize the horizontal coupling across scenarios to $x_{(s)} \leq x$ and $\alpha_{(s)} \leq \alpha$ to reflect the physical reality of resource investment. In production planning, first-stage decisions are sunk-cost capacity acquisitions (e.g., machinery or labor) that act as a ceiling, it is physically illogical to force resource exhaustion if scenario-specific demand is low.

Mathematically, this restricts dual multipliers to be non-negative ($\lambda \in \mathbb{R}_+^{|\mathcal{S}|}$). While numerically more stable, the dual function remains non-differentiable, which is a structural property that triggers the tailing-off effect. In DD, this effect causes the Master Problem’s search to stall, resulting in first-stage candidates that converge too slowly to yield high-quality primal solutions. Our warm-start strategy bypasses this stagnant phase by shifting the search directly into a region of high-quality consensus. The experiments thus serve as a structural proof-of-concept, demonstrating the method's ability to identify high-quality primal solutions even when the dual bound exhibits the slow asymptotic convergence typical of subgradient methods.

\subsection*{Instance generation} \label{sec:instance_generation}
To thoroughly test LOTUS, we require a set of computationally challenging problem instances that adhere to the theoretical requirements and practical motivations outlined in Section \ref{subsec:assum}. Our generation procedure systematically engineers instances with a block-angular structure and structural properties known to produce weak LP relaxations, thereby creating the large integrality gaps that motivate the use of Lagrangian decomposition over standard B\&B solvers. To satisfy the theoretical requirements of compactness and relatively complete recourse, all production capacities and inventory limits are enforced as finite upper bounds, and unmet demand is permitted at a high penalty cost. The generation methodology combines three core sources of complexity:
\begin{enumerate}
    \item \textit{Fixed-cost structures:} High setup costs are included to induce a highly non-convex search space. This directly tests LOTUS's ability to tighten the integrality gap where $Z_{LP}^{\mathcal{S}}$ provides poor guidance.
    \item \textit{Resource scarcity:} A dense matrix of production constraints with tight capacities is designed to make the discovery of a feasible primal bound computationally difficult, necessitating the stable multiplier guidance provided by the LOTUS warm start.
    \item \textit{Scenario heterogeneity:} High-variance stochastic demand scenarios are generated to ensure a complex dual space. This evaluates the framework’s primary practical motivation: identifying a representative active manifold using a computationally cheaper subset.
\end{enumerate}

Based on the preliminary experiments, we fix the warm-start configuration at $10$ iterations on a $30\%$ subset. We use this warm-start configuration in a series of main experiments to perform a head-to-head comparison between LOTUS and a standard DD algorithm \citep{CAROE199937}
across a wide range of problem instances (24 instances each for 5 seeds). All instances are generated to satisfy the theoretical assumptions of compactness and relatively complete recourse, ensuring the stability of the subgradient updates and the feasibility of primal recovery. 

Consistent with industrial profit-maximization objectives, the experiments are performed on the maximization formulation of the production planning problem. In this context, the Lagrangian dual bound serves as a theoretical upper bound, while the feasible solution identified by the projection heuristic constitutes the primal lower bound. This setup allows us to explicitly measure LOTUS's ability to improve the lower bound through accelerated primal convergence, even in scenarios where the upper bound demonstrates characteristic tailing-off effects.

All experiments were carried out on a cluster with a $2.4$GHz CPU with $8$GB RAM, single-threaded, using Gurobi $10.0$ \citep{Gurobi10} in Python $3.9$. The code for generating the production planning problem instances and for running the experiments can be found on GitHub \citep{Cornielje_LOTUS_framework}.

\section{Further Insight Into the Speed Gain} \label{sec:theoretical_insights}
\subsection*{Theoretical Insight of the Speed Gain}\label{subsec:theofound}
The theoretical validity of accelerating the primal convergence can be formalized by analyzing the trade-off between computational cost and information retrieval. Let $t_L(\mathcal{S})$ denote the average time required to perform one subgradient iteration on the scenario set $\mathcal{S}$. The computational advantage of LOTUS compared to DD relies on two structural properties:

\paragraph{1. Cost asymmetry}
The computational cost of a subgradient iteration scales linearly with the number of scenarios. Therefore, iterations performed on the reduced set $\mathcal{S}'$ are significantly cheaper than those on the full set:
\begin{equation}
    t_L(\mathcal{S}') \ll t_L(\mathcal{S}).
\end{equation}
This asymmetry enables LOTUS to perform the necessary ``exploration'' of the dual space at low cost.

\paragraph{2. Active set identification}
The speed gain depends critically on the quality of information transfer between the reduced and full problems. Let $\mathcal{A}_{\mathcal{S}}$ be the set of active non-anticipativity constraints (i.e., indices having non-zero multipliers $\lambda$) in the optimal solution of the full problem. Similarly, let $\mathcal{A}_{\mathcal{S}'}$ be the active set identified by the reduced problem. The net speed gain is maximized when the overlap $|\mathcal{A}_{\mathcal{S}} \cap \mathcal{A}_{\mathcal{S}'}|$ is large. If the reduced set $\mathcal{S}'$ is representative, the warm start $\lambda^{\text{WS}}$ correctly predicts a large portion of $\mathcal{A}_{\mathcal{S}}$. 

Consequently, the number of expensive iterations $N_{\text{full}}$ required to solve the full problem is reduced by the equivalent number of cheap iterations $N_{\text{red}}$ already performed. We can approximate the asymptotic speed gain as:
\begin{equation}
    \text{Gain} \propto \underbrace{N_{\text{red}} \cdot (t_L(\mathcal{S}) - t_L(\mathcal{S}'))}_{\substack{\text{Savings from} \\ \text{cheap exploration}}} \quad- \underbrace{\text{Loss}(\mathcal{S}, \mathcal{S}')}_{\substack{\text{Correction for} \\ \text{misidentified active sets}}}.
\end{equation}
This formulation highlights that LOTUS is most effective when the cost difference $t_L(\mathcal{S}) - t_L(\mathcal{S}')$ is large and the structural overlap between the full problem and the reduced subproblem is high, thereby minimizing the ``Loss'' associated with correcting multipliers that were incorrectly identified as active or inactive in the warm-start phase.

\subsection*{Numerical Results for the Speed Gain}
The efficiency gains for LOTUS compared to DD are quantified in Table \ref{tab:iteration_costs}, which compares the average duration of each iteration. 

\begin{table}[h]
\centering
\caption{Average time (in sec) per subgradient ascent iteration in LOTUS warm-start phase $t_L(\mathcal{S}')$, LOTUS main phase $t_L(\mathcal{S}')$, and DD $t_{DD}(\mathcal{S}')$.}
\label{tab:iteration_costs}
\small 
\begin{tabular}{@{}lrrr|r@{}}
\toprule
& \multicolumn{2}{c}{LOTUS} & \multicolumn{1}{c|}{DD} & Speed gain \\
$|\mathcal{S}|$ & \shortstack{$t_L(\mathcal{S}')$ (sec)} & \shortstack{$t_L(\mathcal{S})$ (sec)} & \shortstack{$t_{DD}(\mathcal{S})$ (sec)} & \shortstack{(\%)}\\ \midrule
$250$  & $17.2$ & $23.5$  & $27.2$  & $13.6$ \\
$500$  & $29.2$ & $47.6$  & $56.3$  & $15.4$ \\
$1000$ & $50.2$ & $89.3$  & $104.3$ & $14.4$ \\
$2000$ & $72.0$ & $162.0$ & $200.0$ & $19.0$ \\ 
\bottomrule
\end{tabular}
\end{table}

At $|\mathcal{S}|=2000$, an iteration on the reduced set costs only $72.0$ sec, while a DD iteration costs $200.0$s. This $2.7$ times cost reduction allows LOTUS to perform its initial dual space exploration in a significantly lower-cost environment. Interestingly, LOTUS demonstrates superior speed even when solving the full scenario set. In Step 5 of LOTUS, the average iteration time is $162.0s$, which is $19\%$ faster than DD with $200.0s$. This suggests that the warm-start multipliers $\lambda^{WS}$ provide a high-quality initial coordination of the solution space that simplifies the solution process for solving the full problem. As the number of scenarios decreases, the speedup of the LOTUS iterations relative to DD iterations decreases slightly. Overall, the iterations in the main phase of LOTUS are $15.8\%$ faster than those in DD.

\begin{table}[h]
\centering
\caption{Allocation of the iteration time budget for LOTUS and DD. The reported times and number of iterations are the averages across all instances (for the full phase of DD, we report the full runtime of $7200$ sec).}
\label{tab:warmstart_iter}
\small 
\begin{tabular}{@{}l|rc|rc|r@{}}
\toprule
 & \multicolumn{2}{c|}{Step 3} & \multicolumn{2}{c|}{Step 5} & \multicolumn{1}{c}{DD} \\
$|\mathcal{S}|$ & \shortstack{time (sec)} & \shortstack{\# iter.} & \shortstack{time (sec)} & \shortstack{\# iter.} & \shortstack{\# iter.}\\ \midrule
$250$  & $172$ & $10$ & $7028$ & $299$ & $265$ \\
$500$  & $292$ & $10$ & $6908$ & $145$ & $128$ \\
$1000$ & $502$ & $10$ & $6698$ & $75$  & $69$  \\
$2000$ & $720$ & $10$ & $6480$ & $40$  & $36$  \\ 
\bottomrule
\end{tabular}
\end{table}

To evaluate the impact of the LOTUS warm-start strategy, we analyze the computational overhead and iteration efficiency as a function of the number of scenarios ($|\mathcal{S}|$). The core of the performance advantage of LOTUS compared to DD lies in the cost asymmetry, see Section \ref{subsec:theofound}, between the reduced scenario set ($\mathcal{S}'$) used in Step 3 (solving the warm start) and the full scenario set ($\mathcal{S}$) used in Step 5 (solving the full problem). Table \ref{tab:warmstart_iter} summarizes how the $120$-minute $(7200s)$ time budget is allocated for both LOTUS and DD. On average, the warm-start phase (Step 3) of LOTUS  consumes $415$ seconds, taking up approximately $5.76\%$ of the total runtime budget. As the number of scenarios increases, the absolute warm-start time grows, reaching $720$ seconds ($10\%$ of the runtime budget) for $|\mathcal{S}|=2000$.
\end{document}